\documentclass[11pt]{article}
\usepackage{fullpage}

\usepackage{epsfig}
\usepackage{amsthm}
\usepackage{amsmath}
\usepackage{amssymb}
\usepackage{algorithm}
\usepackage{algorithmic}
\usepackage{subfig}
\usepackage{graphicx}
\usepackage{color}
\usepackage{epstopdf}
\usepackage[section] {placeins}

\usepackage{hyperref}

\usepackage{natbib}

\usepackage[section] {placeins}

\newcounter{thm_counter}

\newcounter{ass_counter}

\newtheorem{theorem}[thm_counter]{Theorem}
\newtheorem{corollary}[thm_counter]{Corollary}
\newtheorem{assumption}[ass_counter]{Assumption}

\def\rc{\color{black}}

\def\E{\mathbb{E}}
\def\Lmin{\lambda_{\min}}
\def\Lmax{\lambda_{\max}}
\def\Lres{L_{\text{res}}}
\def\eqnok#1{(\ref{#1})}
\def\hogwild{{\sc Hogwild!}}

\newcommand{\RK}{\mbox{\sc RK}}
\newcommand{\AsyRK}{\mbox{\sc AsyRK}}

\newcommand{\CG}{\mbox{\sc CG}}

\newcommand{\R}{\mathbb{R}}
\newcommand{\AsySCD}{\mbox{\sc AsySCD}}

\def\sjwcommentsolved#1{}

\def\jlcommentsolved#1{}

\title{An Asynchronous Parallel Randomized Kaczmarz Algorithm}

\author{Ji Liu, Stephen J. Wright, and Srikrishna Sridhar\\
\{ji-liu, swright, srikris\}@cs.wisc.edu\\
Department of Computer Sciences, University of Wisconsin-Madison}

\begin{document}
\maketitle \thispagestyle{plain} \pagestyle{plain}

\begin{abstract}
We describe an asynchronous parallel variant of the randomized
Kaczmarz ($\RK$) algorithm for solving the linear system $Ax=b$. The
analysis shows linear convergence and indicates that nearly linear
speedup can be expected if the number of processors is bounded by a
constant that depends linearly on the number of rows in $A$ and
inversely on the maximum eigenvalue of $A^TA$.
\end{abstract}

\section{Introduction}
We consider the problem of finding a solution to a consistent linear
system
\begin{equation}
Ax=b, \label{AsyRK:eqn_problem}
\end{equation}
where $A\in \mathbb{R}^{m\times n}$ and $b\in\mathbb{R}^{m}$. We denote the rows and columns of
$A$ by $a_i^T$ and $\bar{a}_j$ respectively, and the elements of $b$ by $b_i$,
$i=1,2,\dotsc,m$, $j=1,2,\dotsc,n$. That is,
 \begin{equation*}
 \begin{aligned}
  A=\left[ \begin{matrix}
      a_1^T \\
      a_2^T \\
      \vdots \\
      a_m^T
   \end{matrix}\right]
  = \left[\begin{matrix}
  \bar{a}_1, \bar{a}_2, \cdots, \bar{a}_n 
  \end{matrix}\right], 
   \quad
   b=\left[ \begin{matrix}
      b_1 \\
      b_2 \\
      \vdots \\
      b_m
   \end{matrix}\right].
      \end{aligned}
\end{equation*}
Besides consistency of $Ax=b$, we assume that throughout that $A$ has
no zero rows.
In fact, we assume (to simplify the analysis) that that the rows of
$A$ are normalized, that is,
\[
\|a_i\|=1,\quad i = 1,2,\dotsc,m,
\]
although we define the algorithm as if normalization had not been
applied. 

We are interested in the casein which $A$ is extremely large and
sparse. The randomized Kaczmarz ($\RK$) is an algorithm for solving
\eqref{AsyRK:eqn_problem} that requires only $O(n)$ storage and has a
linear (geometric) rate of convergence. In some situations, it is even
more efficient than the conjugate gradient ($\CG$) method
\citep{Strohmer09}, which forms the basis of the most popular
iterative algorithms for solving large linear systems. At iteration
$j$, the $\RK$ algorithm randomly selects a row $i(j)
\in\{1,2,\dotsc,m\}$ of the linear system (the probability of choosing
row $i$ is ${\|a_i\|^2}/{\|A\|^2_F}$) and does an orthogonal
projection of the current estimate vector onto the hyperplane {\rc
  $a_{i(j)}^Tx=b_{i(j)}$}:
\begin{equation} \label{eq:onestep}
x_{j+1} = x_{j} - \frac{a_{i(j)}^Tx_j-b_{i(j)}}{\|a_{i(j)}\|^2}a_{i(j)}.
\end{equation}
This update formula can be derived also by applying the basic
stochastic gradient algorithm to the objective $\frac12
\|Ax-b\|^2=\frac12 \sum_{i}(a_i^Tx-b_i)^2$ where $(a_i^Tx-b_i)a_i$ is
the stochastic gradient corresponding to a random choice of index $i$
and $1/\|a_i\|^2$ is the steplength for that gradient estimate. The
expected linear convergence rate of $\RK$ can be proved trivially as
follows \citep{Strohmer09, Needell10, Leventhal08}. Denoting by
$x_j^*$ the projection of iterate $x_j$ onto the solution set of
\eqnok{AsyRK:eqn_problem}, we have from \eqnok{eq:onestep} that 
\begin{align*}
 \|x_{j+1}-x_{j+1}^*\|^2 & \le
\left\| x_j - {1\over \|a_{i(j)}\|^2}a_{i(j)} (a_{i(j)}^T x_j-b_{i(j)}) - x^*_j \right\|^2 \\
&=
\|x_j - x^*_{j} \|^2 - {1\over \|a_{i(j)}\|^2}(a_{i(j)}^Tx_j - b_{i_j} )^2.
\end{align*}
Given the probability $\|a_i\|^2/\|A\|^2_F$ of choosing $i$, we have
by taking expectations that
\begin{align}
\notag
E_{i(j)} \left[ \|x_{j+1}-x_{j+1}^*\|^2 \, | \, x_j \right]
&\le
\|x_j - x_j^* \|^2 - E_{i(j)} \left[{1\over \|a_{i(j)}\|^2} (a_{i(j)}^Tx_j - b_{i(j)} )^2 \right] \\
\notag
&=
\|x_j - x_j^* \|^2 - {1\over \|A\|^2_F}\|Ax_j-b \|^2  \\
\label{eq:linrk}
& \le \left( 1- \frac{\Lmin}{\|A\|^2_F} \right) \|x_j-x_j^* \|^2,
\end{align}
where $\Lmin$ is the smallest nonzero eigenvalue value of $A^TA$.

Recently, asynchronous parallel stochastic algorithms have received
broad attention for solving large convex optimization
problems. \cite{Hogwild11nips} proposed a simple, effective
asynchronous scheme to parallelize the stochastic gradient
algorithm. In this approach, the unknown vector $x$ is stored in
memory locations accessible to all cores of a multicore processor, and
all cores are free to update $x$ in an asynchronous, uncoordinated
fashion. It is assumed that there is a bound $\tau$ on the age of the
updates, that is, no more that $\tau$ updates in total can be occur
between the time at which any processor reads the current $x$
and the time at which it makes its
update. \hogwild~\citep{Hogwild11nips} allows a lock-free
implementation, since the update to a single element of $x$ is an
atomic operation. \cite{Avron13arXiv}, \cite{LiuWright13} and \cite{Sridhar2013nips} applied a
similar asynchronous scheme to stochastic coordinate descent, and have
proved attractive convergence properties.

We apply the same asynchronous parallel technique used in
\hogwild~\citep{Hogwild11nips} to the standard $\RK$ algorithm. The
unknown vector $x$ is stored in a shared location, and all cores
simultaneously run a $\RK$ process, updating $x$ in an asynchronous
fashion.  Although our asynchronous parallel randomized Kaczmarz
algorithm ($\AsyRK$) can be viewed as an application of \hogwild~to
the objective $\frac12 \|Ax-b\|^2$ (with a particular choice of step
length), our analysis shows a linear convergence rate for $\AsyRK$
that outperforms the $1/t$ sublinear convergence rate for \hogwild.

Our analysis also provides an indication of the maximum number of
cores that can be involved in the computation while still yielding
approximately linear speedup. This ``bound'' is expressed in terms of
the number of equations $m$ and the maximal eigenvalue of $A^TA$.

An outline of the remainder of the paper is as follows.  We review
related work in
Section~\ref{sec:relatedwork}. Section~\ref{AsyRK:sec:algorithm}
illustrates details of the $\AsyRK$ algorithm. The convergence rate of
$\AsyRK$ is described in Section~\ref{sec:mainresults}, with proofs
given in Appendix~\ref{AsyRK:sec:proofs}. Some simple experiments
illustrate linear speedup in Section~\ref{AsyRK:sec:exp}. We discuss
extensions to the inconsistent case in Section~\ref{sec:incon}, and
make some concluding observations in
Section~\ref{AsyRK:sec:conclusion}.

\subsection*{Notation and Assumption}

We use the following notation.
\begin{itemize}
\item $\|x\|_0$ denotes the cardinality or ``$\ell_0$ norm'' of the
  vector $x$, that is, the number of nonzero elements in $x$.
\item $\|X\|$ is the spectral norm of the matrix $X$, while $\|X\|_F$
  is the Frobenius norm.
\item $P_{t}$ is the square $n \times n$ matrix of all zeros, except
  for a $1$ in the $(t,t)$ position. 
\item Several quantities characterize the rows and columns of $A$:
  $\theta_i := \|a_{i}\|_0$, $\mu:=\max_i \|a_i\|_0$,
  $\nu:=\max_{j}\|\bar{a}_j\|_0$.
\item $\alpha:=\max_{i,t}\|A\theta_{i}P_ta_i\|$. One can verify that $\alpha\leq \sqrt{\nu}\mu$ and $\alpha\leq \|A\|\mu$.
\item Given $x_j \in \mathbb{R}^n$, $x_{j}^*$ denotes the projection
  of $x_{j}$ onto the solution set of \eqref{AsyRK:eqn_problem}.
\item The support index set of $x$ is defined as $\mbox{\rm supp}(x)$.
\item $\Lmin$ is defined as the minimal \emph{nonzero} eigenvalue value of
  $A^TA$, while $\Lmax$ is defined as the maximal eigenvalue value of
  $A^TA$.
\end{itemize}

We make a few observations about $\Lmax$.  If $A$ is a matrix whose
elements are i.i.d Gaussian random variables from $\mathcal{N}(0,1)$,
then row-normalized, fundamental results in random
matrices~\citep{Vershynin2011arXiv} yield that $\Lmax$ is bounded by
$\left({\sqrt{m}+\sqrt{n}\over \sqrt{n}}\right)^2 \leq O(1+m/n)$ with
high probability. As long as $m/n$ is bounded by a constant, $\Lmax$
is bounded by a constant as well. If $A$ is a sparse matrix, then
\[
\Lmax:=\max_{\|y\|=1}\|A^TAy\| = \max_{\|y\|=1} \|Ay\|^2 \leq \max_{i}
|\{j:~\text{supp}(a_i)\cap \text{supp}(a_j)\neq \emptyset\}| \leq
\mu\nu.
\]

\begin{assumption} \label{ass_e}
Assume that
\begin{itemize}
\item The solution to \eqref{AsyRK:eqn_problem} exists.
\item $A$ is row-normalized, that is, $\|a_i\|=1~\forall
  i\in\{1,2,\cdots,m\}$.
\end{itemize}
\end{assumption}

Note that $\| A\|_F^2=m$ when the rows of $A$ are normalized.

\section{Related Work} \label{sec:relatedwork}

The original Kaczmarz algorithm \citep{Kaczmarz37} used a cyclic
projection procedure to solve consistent linear systems
$Ax=b$. Kaczmarz proved convergence to the unique solution when $A$ is
a square nonsingular matrix. The cyclic ordering of the iterates made
it difficult to obtain iteration-based convergence results, but
\citet{Galantai25} proved a linear convergence rate in terms of
cycles. Since the 1980s, the Kaczmarz algorithm has found an important
application area in Algebraic Reconstruction Techniques (ART) for
image reconstruction; see for example \citet{Herman80, Herman09}. It
is sometimes referred to in this literature as the ``sequential
row-action ART algorithm.''

\citet{Strohmer09} studied the behavior of $\RK$ in the case of a
consistent system $Ax=b$ in which $A$ has full column rank (making the
solution unique).  They proved linear convergence rate for $\RK$ in
expectation.  \cite{Needell10} also assumed full column rank, but
dropped the assumption of consistency, showing that the $\RK$
algorithm converges linearly to a ball of fixed radius centered at the
solution, where the radius is proportional to the distance of $b$ from
the image space of $A$.
\citet{EldarN11} presented a modified version of the randomized
Kaczmarz method which selects the optimal projection from a randomly
chosen set at each iteration. This technique improves the convergence
rate, but requires more computation per iteration.
 \cite{LiuWright13b} proposed an accelerated $\RK$ algorithm that uses
 a Nesterov-type accelerated scheme, improving the linear convergence
 rate constant from $1-\Lmin/m$ 
(corresponding to \eqref{eq:linrk}, after normalization of rows) to
 $1-\sqrt{\Lmin}/ m$.


\citet{Leventhal08} extended the $\RK$ algorithm for consistent linear
equalities $Ax=b$ to the more general setting of consistent linear
inequalities and equalities: $A_Ix\geq b_I$, $A_Ex=b_E$. The basic
idea is quite similar to the $\RK$ algorithm: iteratively update
$x_{k+1}$ by projecting $x_{k}$ onto the hyperplane or half space for
a randomly selected equality or inequality constraint. The linear
convergence rate was proven to be $1-1/{(L^2\|A\|^2_F)}$, where $L$ is
the Hoffman constant \citep{Hoffman52} for the full system. 


\citet{Zouzias12} considered the case of possibly inconsistent
\eqnok{AsyRK:eqn_problem}. They proposed a randomized extended
Kaczmarz algorithm by first projecting $b$ orthogonally onto the image
space of $A$ to obtain $b_\bot$, then orthogonally projecting the
initial point $x_0$ onto the hyperplane $Ax=b_{\bot}$. Essentially,
the $\RK$ algorithm is applied twice. The convergence rate is proven
to be $1-{\Lmin/\|A\|^2_F}$, which is the same as the $\RK$ algorithm
for consistent linear systems. This method can be considered as a
randomized variant of the extended Kaczmarz method proposed by
\cite{Popa99}.

Among synchronous parallel methods, \citet{Cegogo01} proposed a
parallel component averaging method to
solve~\eqref{AsyRK:eqn_problem}. This approach parallel-projects the
current $x$ onto all (or multiple) hyperplanes, then applies an
averaging scheme to the projections to obtain the next iterate. This
method is essentially a gradient descent method for solving ${1\over
  2}\|Ax-b\|^2$, so is able to handle inconsistent problems.  This
paper also notes \citep[Section~5.1]{Cegogo01} that for sparse
problems, parallelism can be obtained by simultaneously projecting the
current iterate onto a set of mutually orthogonal hyperplanes,
obtained by considering equations whose nonzero components appear in
disjoint locations. When obtained forom image reconstruction problems,
such sets of equations can be obtained by considering parallel rays
that are sufficiently far apart so as to pass through disjoint sets of
pixels. This type of parallelism has small granularity, and the amount
of communication required between processors may make it unattractive
in practice.

A related approach is Block-Cimmino (or Block-AMS) algorithm
\citep{Aharoni89}, which can be considered as the block version of
Kaczmarz algorithm. Other variants of Block-Cimmino algorithm are
described in \citep{ElfvingNikazad09, Nikazad08}.

Another synchronous parallel approach (for general convex
optimization) due to \citet{Ferris94} distributes variables among
multiple processors and optimizes concurrently over each subset.  A
synchronization step searches the affine hull formed by the current
iterate and the partial optima found by each processor.

In discussing asynchronous parallel methods, we make a distinction
according to whether it is assumed that the reading of $x$ by each
processor is ``consistent'' or not. The term ``consistent'' in this
context means that the $x$ used by each processor to evaluate its
update is an iterate that actually existed at some point in time,
whose components were not changed repeatedly by other processors
during reading (yielding a hybrid of two or more iterates).
\cite{Bertsekas89} introduced an asynchronous parallel implementation
for general fixed point problems $x = q(x)$ over a separable convex
closed feasible region. The optimization problem of minimizing $f$
over a closed convex set $\Omega$ can be formulated as a fixed-point
problem by defining $q(x) := \mathcal{P}_\Omega[(I - \alpha \nabla
  f)(x)]$, where $\mathcal{P}_\Omega$ denotes Euclidean projection
onto $\Omega$. The vector $x$ is stored in memory accessible to all
cores, and the cores update the value of $x$ without locking or
coordination. Inconsistent reading of $x$ is allowed. Linear
convergence is established --- using admirably straightforward
analysis --- provided that $\nabla^2 f(x)$ satisfies a diagonal
dominance condition, guaranteeing that the iteration $x=q(x)$ is a
maximum norm contraction mapping for sufficient small $\alpha$. We
note, however, that this condition is even stronger than strong
convexity.

\citet{Elsner90} proposed an asynchronous parallel $\RK$ algorithm,
again for a situation in which all processors have access to $x$
stored in commonly accessible memory. Each processor iteratively runs
the following procedures, where $x$ denotes a globally shared version
of the variable vector and $x'$ and $x''$ denote copies stored locally
on each processor: (a) read the current global $x$ into the local
$x'$; (b) write the convex combination of the local variables $x'$ and
$x''$ into $x''$, and also into the shared memory as a new $x$; (c)
project the local $x''$ onto the selected hyperplane to get a new
$x''$.  Note that the algorithm requires locking the shared memory in
step (b) because it does not allow two processors to access shared
memory at the same time. Our computational experiences with related
algorithms (e.g., $\AsySCD$ \citep{LiuWright13} and
\hogwild~\citep{Hogwild11nips}) indicate that memory locking of this
type degrades computational performance seriously. Moreover, the
convergence analysis establishes convergence, but does not prove a
linear convergence rate.

\hogwild~\citep{Hogwild11nips} is a lock-free, asynchronous parallel
version of the stochastic gradient method. All processors share the
same memory storing $x$ and update it simultaneously. Unlike
\cite{Bertsekas89}, inconsistent reads of $x$ are not permitted by the
analysis.
When the updates satisfy a certain sparsity property, the convergence
of \hogwild~ approximately matches the $1/t$ rate of serial stochastic
gradient, as described and analyzed by \citep{Nemirovski09}. Recent
work by \cite{Avron13arXiv} concerned an asynchronous linear solver
for $Ax=b$ (for symmetric positive definite $A$) using the same
asynchronous scheme as \hogwild, proving a linear convergence rate.

\cite{LiuWright13} followed the model of \hogwild~to propose an
asynchronous parallel stochastic coordinate descent ($\AsySCD$)
algorithm and proved sublinear ($1/t$) convergence on general convex
functions and a linear convergence rate on functions that satisfy an
essential strong convexity property.
\cite{Richtarik12arXiv} proposed a parallel coordinate descent method
for minimization of a composite convex objective with separable
nonsmooth part. Their method is a synchronous parallel approach
\sjwcommentsolved{Is it true that the {\em implemenation} of this
  method is asynchronous? We should say this, if so. (Yes. I checked
  their paper. It analyzed a synchronous method but implemented it in
  an asynchronous way.)} (in contrast to $\AsySCD$, which is
asynchronous), but it is implemented in an asynchronous
fashion. Another distinction between the two approaches is found in
the convexity assumptions, which are slightly weaker in
\cite{LiuWright13}.





\section{Algorithm} \label{AsyRK:sec:algorithm}

Each thread in our $\AsyRK$ algorithm performs the following simple
steps: (1) Choose\jlcommentsolved{``Choose'' or ``choose''?} an index
$i$ randomly from $\{1,2,\dotsc,m\}$; (2) read the components of $x$
that correspond to the nonzeros in $a_i$ from shared memory; (3)
calculate $a_i^Tx-b_i$; (4) select $t \in \text{supp}(a_i)$; (5)
update component $t$ of $x$ in the shared memory by a multiple of
$(a_i)_t (a_i^Tx-b_i)$.  In principle, no memory locking takes place
during either read or write, but we assume that the reads are
``consistent,'' according to the discussion above. (We note that
inconsistent reading is expected to be rather rare in the case of
sparse $A$, because only those elements of $x$ that correspond to
nonzero locations in $a_i$ need to be read, and inconsistency possibly occurs
only when this subset of elements is updated at least twice by other
processors while it is being read.) 
The update to component $t$ of $x$ can be implemented as a unitary
operation, requiring no memory locking.

Algorithm~\ref{alg_rka} gives a global, aggregated view of this
multithreaded process. An iteration counter $j$ is incremented each
time $x$ is updated by a thread. We use $k(j)$ to denote the iterate
at which $x$ was read by the thread that updated $x_j$ to
$x_{j+1}$. (We always have $k(j) \le j$, and strict inequality holds
when other threads have updated $x$ between the
time\jlcommentsolved{``times'' or ``time''?} it is read and the time
the update is performed by this thread.)  The index $i(j) \in
\{1,2,\dotsc,m\}$ denotes the row that was selected by the thread that
updated $x_j$ to $x_{j+1}$. The index $t(j) \in \text{supp}
(a_{i(j)})$ denotes the component of $x$ that is chosen (randomly) to
be updated at iteration $j$.  We assume that the delay between reading
and update for each thread is not too long, that is,
\begin{equation} \label{eq:kj}
k(j) \ge j-\tau,
\end{equation}
for some integer $\tau \ge 1$.  $\tau$ can be assumed to be similar to
the number of processors that are involved in the computation. Note
that the step depends on $\theta_{i(j)}$ (the cardinality of the
chosen row) and a parameter $\gamma$ which is critical to the analysis
of the following sections.





\begin{algorithm}[htp!]               
\caption{Asynchronous Randomized Kaczmarz Algorithm
  $x_{K+1}=\AsyRK(A,b,x_0, \gamma, K)$} 
\label{alg_rka}                           
\begin{algorithmic}[1]                    
\STATE Given $A \in \R^{m \times n}$ with normalized rows $a_i$,
$i=1,2,\dotsc,m$, $b \in \R^m$;
\STATE Initialize $j\leftarrow 0$;
\WHILE{$j \leq K$} \STATE Choose $i(j)$
from $\{1,2,\dotsc,m\}$ with equal probability;
\STATE Choose $t(j)$ from ${\rm supp}(a_{i(j)})$ with equal probability; \label{eqn_alg_rka_2} 
\STATE Update
$x_{j+1}\leftarrow x_j-\gamma
\theta_{i(j)}P_{t(j)}a_{i(j)}(a_{i(j)}^Tx_{k(j)}-b_{i(j)})$;
\label{eqn_alg_rka_3} \STATE $j\leftarrow j+1$; \ENDWHILE
\end{algorithmic}
\end{algorithm}

\section{Main Results} \label{sec:mainresults}

This section presents the convergence analysis for $\AsyRK$. The key
issue for $\AsyRK$ is to choose an appropriate steplength parameter
$\gamma$.  At an intuitive level, we would like $\gamma$ to be large
enough to make significant progress in the approximate gradient
direction. On the other hand, we want to keep it small enough that the
approximate gradient information computed at the earlier iterate
$k(j)$ is still relevant when the time comes to do the update at
iteration $j$. That is, the difference between $x_{k(j)}$ and $x_j$
should not be too large. Along these lines, we require the ratios of
expected residuals at any two successive iterations to be bounded
above and below, as follows:
\[
\rho^{-1} \leq \frac{\E\|Ax_{j+1}-b\|^2}{\E\|Ax_{j}-b\|^2} \leq \rho,
\]
where $\rho$ is a user defined parameter, usually set to be slightly
larger than $1$. The steplength $\gamma$ depends strongly on $\rho$.

We state a result about convergence of $\AsyRK$ in
Algorithm~\ref{alg_rka}.
\begin{theorem} \label{AsyRK:thm_2}
Assume that Assumption~\ref{ass_e} is satisfied. Let $\rho$ be any
number greater than 1 and define the quantity $\psi$ as follows:
\begin{equation}
\psi:=\mu + \frac{2\Lmax\tau\rho^{\tau}}{m}.
\label{eqn_thm_psi}
\end{equation}
Suppose the steplength parameter $\gamma>0$ in Algorithm~\ref{alg_rka}
satisfies the following three bounds:
\begin{equation}
\gamma\leq {1\over \psi},\quad\gamma\leq \frac{m(\rho-1)}{2\Lmax\rho^{\tau+1}},\quad\gamma\leq m\sqrt{\frac{(\rho-1)}{\rho^{\tau}(m\alpha^2+\Lmax^2\tau\rho^{\tau})}}.\label{eqn_thm2_gamma}
\end{equation}
Then we have for any $j\geq 0$ that
\begin{equation}
\rho^{-1}\E(\|Ax_j-b\|^2)\leq \E(\|Ax_{j+1}-b\|^2)
\leq \rho\E(\|Ax_j-b\|^2) \label{AsyRK_eqn_thm2_1}
\end{equation}
and
\begin{equation}
\E(\|x_{j}-x^*_{j}\|^2) \leq
\left(1-{\Lmin \gamma\over m}(2-\gamma\psi)\right)^j\|x_0-x_0^*\|^2.
\label{AsyRK:eqn_thm2_2}
\end{equation}
\end{theorem}
This theorem indicates a linear rate of convergence, outperforming the
sublinear ``$1/j$'' convergence rate for the asynchronous stochastic
gradient method \hogwild. The key reason for this improvement is that
because $a_i^Tx^*-b_i=0$ for all $i$, the stochastic gradient {\em
  estimates} all approach zero as $x$ approaches $x^*$, a property
that does not hold for general stochastic gradient algorithms.
\jlcommentsolved{I may not agree with this argument. Even given the
  condition (6), \hogwild cannot give linear convergence rate, because
  a special case (the standard stochastic gradient method) does not
  guarantee the linear convergence rate for strongly convex
  function. The true reason is that in our problem we have the
  stochastic gradient $a_i^Tx^* - b_i =0$ for any optimal solution
  $x^*$. In comparison, one does not have this property in for a
  general problem using stochastic gradient method. However, this
  intuition is buried in the proof. SJW: OK, I modified the
  discussion.}

Note that the upper bound on steplength parameter $\gamma$ decreases
as the bound $\tau$ on the age of the iterates increases. This
dependency allows us to figure out how many threads can be executed in
parallel without significantly degrading the convergence behavior.

This following corollary proposes an interesting particular choice for
the parameters for which the convergence expressions become more
comprehensible. The result requires a condition on the delay bound
$\tau$ in terms of $m$ and $\Lmax$.
\begin{corollary} \label{cor}
Suppose that Assumption~\ref{ass_e} is satisfied and that
\begin{equation}
\frac{2e\Lmax(\tau+1)}{m} \leq 1.
\label{eqn_cor_ass}
\end{equation}
Then if we choose
\begin{equation}
\rho=1+\frac{2e\Lmax(\tau+1)}{m}
\label{eqn_cor_rho}
\end{equation}
and set $\gamma = 1/\psi$, where $\psi$ is defined as
in~\eqref{eqn_thm_psi}, we have that
\begin{equation}
\E(\|x_{j}-x^*_{j}\|^2) \leq \left(1 - \frac{\Lmin}{m(\mu+1)}\right)^j \|x_0-x_0^*\|^2.
\label{eqn_cor_estbound}
\end{equation}
\end{corollary}
Over a span of $m$ iterations, \eqref{eqn_cor_estbound} implies a
decrease factor of approximately $1-{\Lmin}/{(\mu+1)}$. This rate
estimate indicates that for a delay $\tau$ (and hence a number of
processors) in the range implied by~\eqref{eqn_cor_ass}, the number of
iterations required for convergence is not affected much by the delay,
so we can expect an almost linear speedup from the multicore
implementation in this regime.

We conclude this section with a high-probability estimate for convergence of $\{\|x_j-x_j^*\|^2\}_{j=1,2,\dots}$.
\begin{theorem} \label{thm_3}
Suppose that the assumptions of Corollary~\ref{cor} hold, and that $\rho$ and $\psi$ are defined as there. For $\epsilon>0$ and $\eta\in (0, 1)$, if
\begin{align}
j \geq \frac{m(\mu+1)}{\Lmin} \left|\log{\|x_0-x_0^*\|^2\over \eta\epsilon}\right|
\label{eqn_thm3_j}
\end{align}
we have that
\begin{align}
\mathbb{P}(\|x_j-x_j^*\|^2\leq \epsilon) \geq 1-\eta.
\end{align}
\end{theorem}

The proofs of all results in this section appear in
Section~\ref{AsyRK:sec:proofs}.

\section{Comparison} \label{sec:cmp}

\begin{table} [h]
\caption{Comparison among $\RK$, $\AsySCD$, and $\AsyRK$. The quantity
  $\delta$ is the fraction of nonzero entries in $A$, while $\Lres$ is
  the maximal row norm of $A^TA$ and $L_{\max}$ is the maximal
  diagonal entry of $A^TA$. (We assume that the nonzeros are roughly
  evenly distributed in $A$, so that $\mu$ is a modest multiple of
  $\delta n$.) The first row shows the number of operations required
  per iteration. The linear convergence rate (in the sense of
  iterations) is shown in the second row. The third row shows the
  maximum number of cores for which linear speedup is available. The
  fourth row combines the preceding rows to obtain the convergence
  rate in the sense of running time, when the method is run on the
  ``maximal'' number of processors.}
\begin{center}
\begin{tabular} {|l|ccc|}
\hline
algorithms & $\RK$ & $\AsySCD$ & $\AsyRK$
\\ \hline 
$\#$ operation per iteration & $O(\delta n)$ & $\min\{O(\delta^2mn),~O(n)\}$ & $O(\delta n)$
\\ 
rate (iteration) & $1-\frac{\Lmin}{m}$ & $1- \frac{\Lmin}{2n L_{\max}}$ & $1-\frac{\Lmin}{m(\mu+1)}$
\\ 
$\#$ processors & $1$ & $O\left(\frac{\sqrt{n}L_{\max}}{\Lres}\right)$ & $O\left(\frac{m}{\Lmax}\right)$
\\ 
rate (running time) & $1-O\left(\frac{\Lmin}{\delta mn}\right)$ & $1 - O\left(\frac{\Lmin}{n^{1.5}\Lres\min\{\delta^2m,~1\}}\right)$ & $1- O\left(\frac{\Lmin}{\delta^2 n^2\Lmax}\right)$
\\ \hline
\end{tabular} \label{AsyRK:tab_cmp}
\end{center}
\end{table}

This section compares the theoretical performance of $\RK$, $\AsySCD$
(applied to minimization of $\frac12 \|Ax-b\|^2$) and $\AsyRK$.  In
Table~\ref{AsyRK:tab_cmp}, we show the complexities (per iteration)
and convergence rates (with respect to number of iterations) of three
algorithms in the first and second rows. The third row gives the
maximal possibly number of processors to parallelize three algorithms
respectively. The last row computes the convergence rate in term of
the operation using the possibly maximal number of processors, which
can be roughly understood as the running time comparison.  In
reporting statistics for $\AsySCD$, we consider two alternative
implementations: (1) randomly choose a coordinate $i$ and compute it
by $(a_iA)x$, which needs $O(\delta^2mn)$ operations per iteration;
and (2) compute $A^TA:=Q$ offline and randomly choose an coordinate
$i$ to compute it by $Q_{i.}x$, which needs $O(n)$ operations per
iteration. We report the complexity per iteration of $\AsySCD$ as the
minimum of these two estimates.

We perform a comparison of convergence behavior of these three
algorithms on a Gaussian random matrix $A$ with i.i.d. elements
generated from $\mathcal{N}(0, 1/n)$. All rows of $A$ have norm
approximately $1$, and $\Lmax$ is approximately $1+m/n$.\jlcommentsolved{seems to be $1+m/n$} For these
values, the convergence rates per iteration of $\AsySCD$ and $\RK$ are
similar. By comparison, the convergence rate of $\AsyRK$ seems worse
than $\RK$ and $\AsyRK$ by a factor $(\mu+1)$. This is because
$\AsyRK$ only updates a single coordinate rather than all coordinates
corresponding to the nonzero elements in the stochastic gradient. If
we modify Algorithm~\ref{alg_rka} to updated {\em all} components in
$\mbox{\rm supp}(a_{i(j)})$ (rather than just the component $t(j)$),
the convergence rate for $\AsyRK$ becomes quite similar to the other
two methods, without an appreciable increase in cost.


Next we compare the parallel implementations.  From the last row of
Table~\ref{AsyRK:tab_cmp}, we see that $\AsyRK$ improves the rate of
$\RK$ if $\delta mn \gg \delta^2 n^2\Lmax$, or equivalently $m\gg \delta
n\Lmax$. Assuming the Gaussian ensemble for $A$ and that $m$ and $n$
are comparable (so that $\Lmax = O(1)$), there is a potential factor
of improvement in runtime of $O(1/\delta)$ for $\AsyRK$ over $\RK$. To
compare $\AsyRK$ and $\AsySCD$, we note that $\Lmax=O(\Lres)$ under
the same scenario for $m$, $n$, and $A$. Comparing the rates (running
time) in the last row of Table~\ref{AsyRK:tab_cmp}, we find that when
the mild condition $\delta < O(n^{-1/4})$ holds, $\AsyRK$ converges
much faster than $\AsySCD$.  Overall, $\AsyRK$ has a clear advantage
in complexity when applied to sparse problems.

\section{Experiments} \label{AsyRK:sec:exp}

We illustrate the behavior of $\AsyRK$ on sparse synthetic data. Our
chief interest is the efficiency of multicore implementations (one
thread per core), compared to a single-thread implementation.

To construct a sparse matrix $A\in\R^{m\times n}$, given dimensions
$m$ and $n$ and sparsity ratio $\delta$, we select $\delta mn$ entries
of $A$ at random to be nonzero and $\mathcal{N}(0,1)$ normally
distributed, and set the rest to zero. Finally, the rows of $A$ are
normalized.


Our experiments run on $1$ to $10$ threads on an Intel Xeon machine,
with all threads sharing a single memory socket. Our implementations
deviate modestly from the version of $\AsyRK$ analyzed here.  First,
$A$ is partitioned into slices (row submatrices) of equal size, and
each thread is assigned one slice. Each thread then selects the rows
in its slice to update in order, with the order being reshuffled after
each scan. This scheme essentially changes from sampling with
replacement (as analyzed) to sampling without replacement, which has
empirically better performance. (The same advantage is noted in
implementations of \hogwild.) The second deviation from the analyzed
version is that {\em all} coordinates corresponding to nonzeros in the
selected row $a_{i(j)}$ are updated, not just the $t(j)$ component.
This scheme makes a single thread behave like $|a_{i(j)}|$ threads,
thus implicitly increasing the number of cores involved in the
computation. Note that this variant represents the obvious extension
of randomized $\RK$. In fact, when implemented on a single thread, it
is precisely the usual randomized $\RK$ scheme. \jlcommentsolved{Shall we
  emphasize that using single thread in $\AsyRK$ is exactly the same
  as the standard $\RK$? SJW: OK, modified this discussion.}


\begin{figure}
  \centering
 \includegraphics[width=0.4\textwidth]{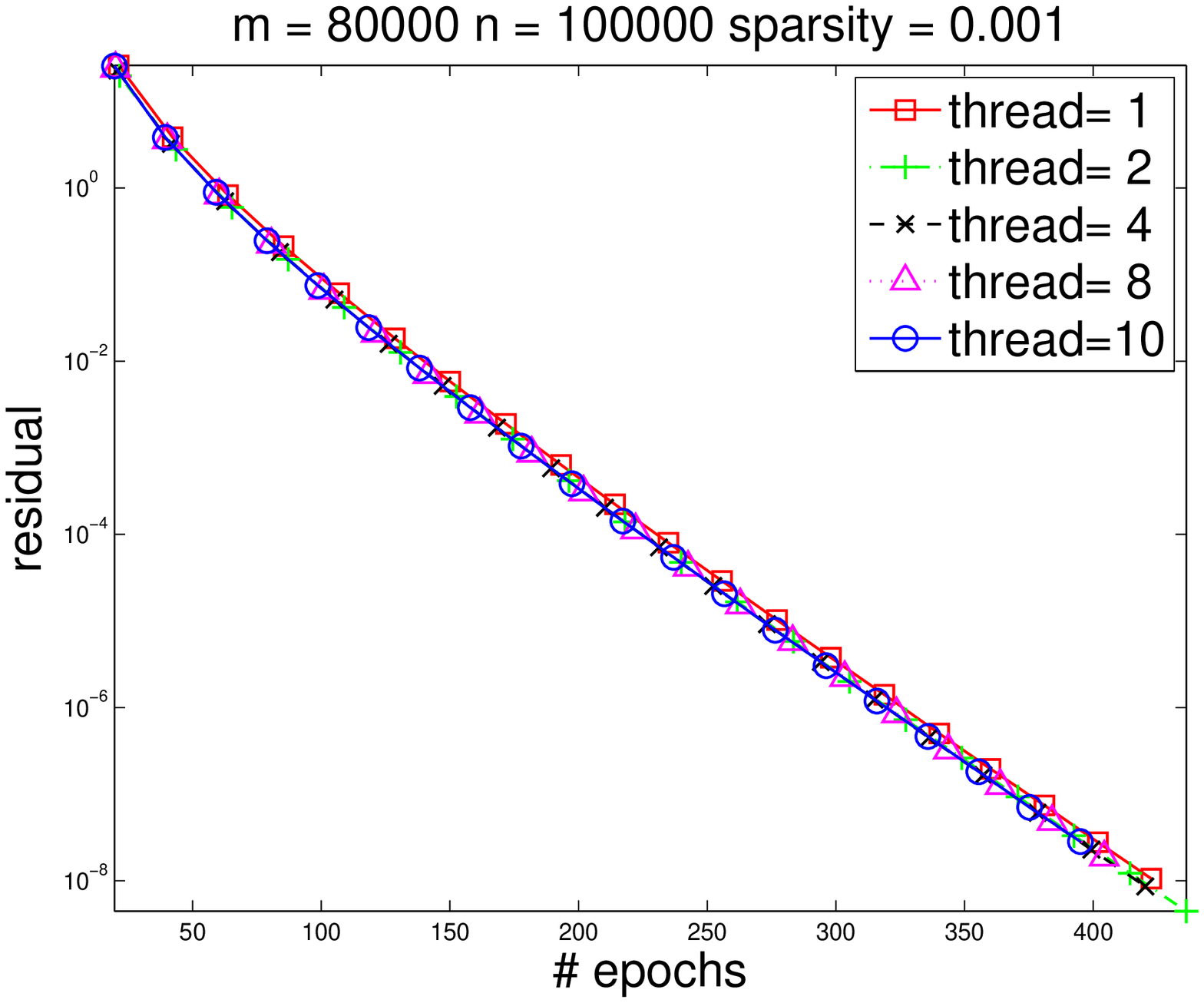} \;\;
 \includegraphics[width=0.4\textwidth]{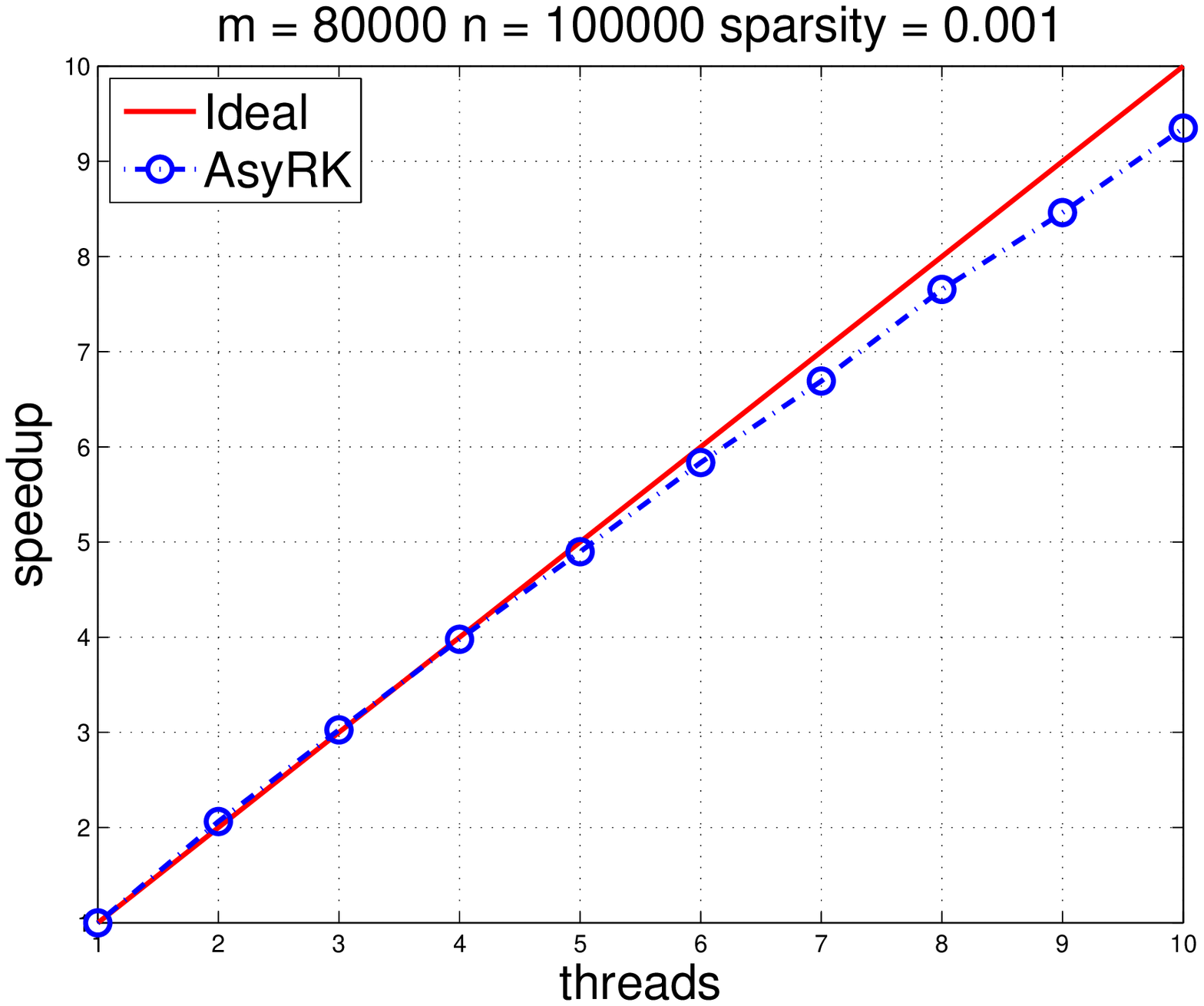}
\caption{The left figure shows one line for each number of threads
  ($=$cores), plotting squared residual $\|Ax-b\|^2$ against epochs. The right
  figure shows the speedup over different numbers of cores.}
\label{fig:ARK_1}
\end{figure}

\begin{figure}
  \centering
 \includegraphics[width=0.4\textwidth]{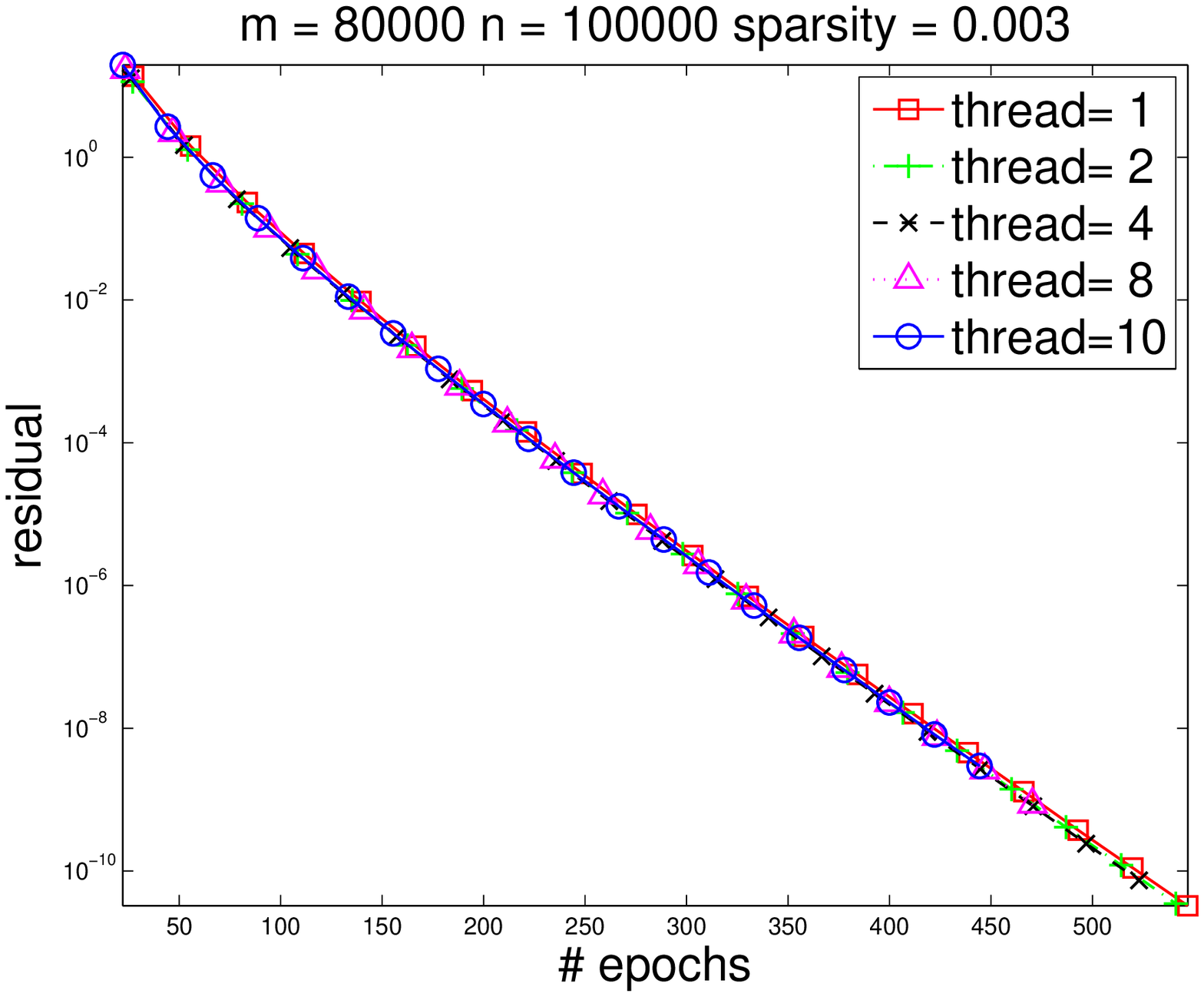} \;\;
 \includegraphics[width=0.4\textwidth]{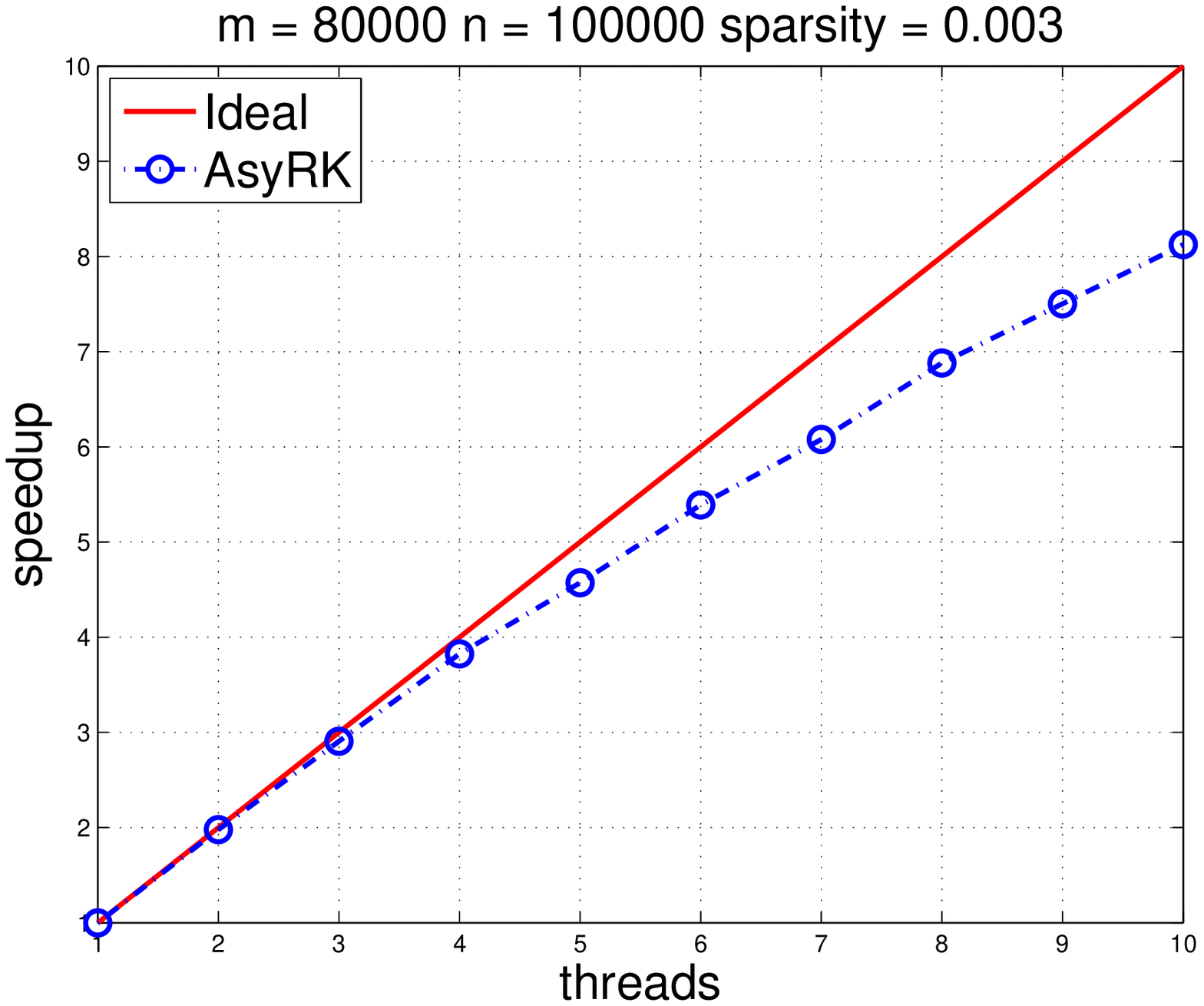}
\caption{The left graph shows one line for each number of threads
  ($=$cores), plotting residual $\|Ax-b\|^2$ against epochs. The
  right graph shows the speedup over different numbers of cores.}
\label{fig:ARK_2}
\end{figure}

For the plots in Figures~\ref{fig:ARK_1} and \ref{fig:ARK_2}, we
choose $m=80000$ and $n=100000$, with $\delta = 0.001$, and set the steplength $\gamma$ as $1$ in
Figure~\ref{fig:ARK_1} and $\delta=0.003$ in Figure~\ref{fig:ARK_2}.
The left-hand graph in each figure indicates the number of threads /
cores and plots residual (defined as
$\|Ax-b\|^2$)\sjwcommentsolved{This is not really a conventional
  definition of ``residual.'' What about $\|Ax-b\|^2$? Or $\| Ax-b\|$?
  (I can agree with you. I discussed this with Krishna before. He used
  the same measure in his $\AsyRK$ implementation as well. The main
  reason why we use the gradient norm to evaluate the residual is that
  the $\AsySCD$ implementation uses this quantity as the residual
  since it only has the information of $Q$. To make all methods
  comparable, we chose to use this strange measurement.) {\bf SJW: I
    suggest for the two figures to use a more conventional definition
    i.e. $\|Ax-b\|$. You'll need a different convergence tolerance but
    the figures should look roughly the same. For the table
    (comparison with AsySCD) you can still use the previous silly
    measure.} (JL: I replaced two graphs using the residual you suggested and also changed captions for figures 1 and 2 and table 2 accordingly. Could you check with it? BTW, I used the Matlab in my desktop to re-draw these two graphs. It looks different from our previous graphs for the Matlab inconsistency issue. I will reproduce them using the laptop you give it to me later on.)} 
vs epoch count, where one epoch is equivalent to $n$
iterations. Note that the curves tend to merge, indicating that the
workload required for $\AsyRK$ is almost independent of the number of
cores. This observation validates our result in Corollary~\ref{cor},
which indicates that provided it is below a certain threshold, the
value of $\tau$ does not affect convergence rate.
The right-hand graph in each figure shows speedup over different
numbers of cores. Near-linear speedup is observed for $\delta=0.001$
(Figure~\ref{fig:ARK_1}), while for $\delta=0.003$ there is a dropoff
for larger numbers of cores (Figure~\ref{fig:ARK_2}).  This can
perhaps be explained by the difference between our implementation from
the analyzed version, in that the nonzeros in the full row $a_{i(j)}$
are updated rather than just a single element. The effect of this
policy can be incorporated into the analysis roughly by increasing the
value of the maximum delay parameter $\tau$. In this case, a matrix
that is three times more dense could be modeled by a value of $\tau$
that is three times larger. The effect may be to raise $\tau$ above
the threshold for which linear speedup can be expected, thus
explaining the (graceful) degradation in speedup for larger numbers of
cores in Figure~\ref{fig:ARK_2}.


\sjwcommentsolved{I don't understand this stuff about ``3 times more
  cores.'' Need to be much clearer here about what is going on. (Sorry
  for this confusion. Here, I want to explain why the speedup curve
  for $\delta = 0.003$ is worse than $\delta=0.001$. Our theory only
  suggests the dependence between speedup and two quantities: $\Lmax$
  and $\tau$ (the number of cores). As we analyzed before, the value
  of $\Lmax$ for $\delta=0.001$ and $\delta=0.003$ are close. Also in
  both cases, we use the same number of physical cores, i.e., 10
  cores. It means that the sparsity ratio should be affect the speedup
  at all. Here I was trying to explain that the number of cores used
  in $\delta=0.003$ is more than $\delta=0.001$ because of our
  implementation trick. In fact, ``$\delta=0.003$'' uses $0.003* n*10$
  cores, while ``$\delta=0.001$'' uses $0.001* n*10$ cores .  That is
  why I said ``$\delta=0.003$ involves 3 times more cores than
  $\delta=0.001$''.) {\bf SJW: OK. I wrote some new sentences to
    explain this. Let me know what you think.} (Thanks. That is good.)}

Next, we compare $\AsyRK$ to $\AsySCD$ \citep{LiuWright13} on sparse
synthetic data sets, on $10$ cores (single socket) of the Intel
Xeon. Various values of $m$, $n$, and $\delta$ are chosen for
comparison in Table~\ref{tab:ARK_1}. A similar number of epochs is
required by both algorithms, reflecting the similarity of their
theoretical convergence rates; see Section~\ref{sec:cmp}. However,
$\AsyRK$ is one order of magnitude faster than $\AsySCD$ to achieve
the same accuracy. The main reason is that, as we showed in
Table~\ref{AsyRK:tab_cmp}, the per-iteration complexity of $\AsySCD$
is much higher than $\AsyRK$, for these values of the parameters.

\begin{table} 
  \centering
\caption{Comparison of running time and epochs between $\AsySCD$ and
  $\AsyRK$ on $10$ cores. We report their running time and number of
  epochs required to attain a residual of $10^{-5}$, where the residual is defined by $\|A^T(Ax-b)\|^2$ for purposes of comparison.}
\begin{tabular}{|rrr|c|rr|rr|}
\hline 
\multicolumn{3}{|c|}{synthetic data} 	&	size (MB)	& \multicolumn{2}{c|}{running time (sec)}		& \multicolumn{2}{c|}{epochs}
 \\
\hline
$m$	& $n$ & $\delta$ &		& $\AsySCD$	& $\AsyRK$ 	& $\AsySCD$	& $\AsyRK$ 
\\
\hline
80000 & 100000 & $0.0005$ 	& 43 &	39.		& 3.6		& 199	& 195
\\
80000 & 100000 & $0.001$ 	& 84 & 	170.	 	& 7.6		& 267	& 284
\\
 80000 & 100000 & $0.003$ 	& 244 &	1279.	& 18.4		& 275	& 232
\\
500000  & 1000000  & $0.00005$	& 282 &	54.		& 5.8		& 19	& 19
\\
500000 & 1000000 & $0.0001$	& 550 &	198.		& 10.4	& 24	& 30
\\
 500000  & 1000000 & $0.0002$	& 1086 &	734.		& 15.0	& 29	& 31
\\
\hline
\end{tabular}
\label{tab:ARK_1}
\end{table}

\section{Extension to Inconsistent Systems}
\label{sec:incon}

Although this paper assumes that the linear system is consistent, we
can extend the algorithm described above to find the least-squares
solution of inconsistent linear systems.

The minimizer of the least-squares objective $\|Ax - b\|^2$ is
equivalent to the linear system $A^TAx = A^Tb$, which can be stated as
the following square, consistent system of linear equations:
\begin{equation} \label{eq:inc1}
Ax - \zeta y = 0,\quad \phi A^Ty = \phi A^Tb,
\end{equation}
for any positive values of $\zeta$ and $\phi$.  Similar reformulations
have appeared previously in the literature; see \citet{Eggermont81},
for example. Here we are mainly interested in the optimal values for
$\zeta$ and $\phi$. We can choose $\zeta$ and $\phi$ to maximize the
critical quantity in the analysis of Algorithm~\ref{alg_rka}, which is
the ratio of the minimum nonzero eigenvalue of $A^TA$ to its squared
Frobenius norm. (In Theorem~\ref{AsyRK:thm_2}, this ratio appears as
$\Lmin/m$, because of the normalization of the rows of $A$.) To show
how this quantity depends on $\zeta$ and $\phi$, we note first that
the coefficient matrix in \eqref{eq:inc1} is
 \[
 \tilde{A}=\begin{bmatrix}
0 & \phi A^T \\
A & -\zeta I
\end{bmatrix}.
\]
Denoting the nonzero singular values of $A$ by $\sigma_1\geq \sigma_2
\geq \cdots \geq \sigma_r>0$ and the full SVD of $A$ by $A=U\Sigma
V^T$ where $U\in \R^{m\times m}$, $V\in \R^{n\times n}$, and $\Sigma
\in \R^{m\times n}$, we can decompose $\tilde{A}$ as follows:
\[
\begin{bmatrix}
V & 0 \\
0 & U
\end{bmatrix}
\begin{bmatrix}
0 & \phi\Sigma^T \\
\Sigma & -\zeta I
\end{bmatrix}
\begin{bmatrix}
V^T & 0 \\
0 & U^T
\end{bmatrix}.
\]
The singular values of $\tilde{A}$ are identical to the singular
values of the center matrix, which can be written after symmetric
permutation as follows:
\[
\begin{bmatrix}
0 & \phi\sigma_r \\
\sigma_r & -\zeta \\
& &\ddots  \\ 
& & & 0 & \phi\sigma_1  \\
& & & \sigma_1 & -\zeta \\
& & & & & -\zeta I \\
& & & & & & \mathbf{0}
\end{bmatrix}.
\]
Hence, the minimal nonzero singular value of $\tilde{A}$ is
\[
\min \left\{ \zeta,~-{\zeta\over 2} + {1\over 2}\sqrt{\zeta^2 + 4\phi \sigma^2_r}\right\}.
\]
Noting that $\|\tilde{A}\|_F^2 = (1+\phi^2)\|A\|_F^2 + m\zeta^2$, we
optimize the critical ratio by finding $\zeta$ and $\phi$ that maximize
the quantity
\begin{equation}
\frac{\min \left\{ \zeta,~-{\zeta\over 2} + {1\over 2}\sqrt{\zeta^2 + 4\phi \sigma^2_r}\right\}^2} {(1+\phi^2)\|A\|_F^2 + m\zeta^2} 
= \min\left\{\frac{\zeta^2} {(1+\phi^2)\|A\|_F^2 + m\zeta^2},~\frac{\left\{-{\zeta\over 2} + {1\over 2}\sqrt{\zeta^2 + 4\phi \sigma^2_r}\right\}^2}{(1+\phi^2)\|A\|_F^2 + m\zeta^2} \right\}.
\label{eqn_inconsistent_1}
\end{equation}
For fixed $\phi$, the first term is monotonically increasing with
respect to $\zeta$ while the second term is monotonically decreasing
with respect to $\zeta$. We can thus express the optimal $\zeta$ value
as $\zeta^*=\sigma_r\sqrt{\phi/2}$, which is the value for which these
two terms are equal. By substituting this value into
\eqref{eqn_inconsistent_1}, we obtain
\begin{equation}
\frac{\sigma_r^2 \phi}{2 (1+\phi^2)\|A\|_F^2 + m\sigma_r^2 \phi} =
\frac{\sigma_r^2}{2 (1/\phi+\phi)\|A\|_F^2 + m\sigma_r^2}.
\label{eqn_inconsistent_2}
\end{equation}
It is clear from the last expression that the optimal value for $\phi$
is $\phi^* = 1$, giving the following maximal value for
\eqref{eqn_inconsistent_2}:
\[
\frac{{\sigma_r^2}} {4\|A\|_F^2 + {m\sigma_r^2}}.
\] 

We conclude that by normalizing the rows of $A$, estimating its
minimum singular value $\sigma_r$, and setting $\phi=1$ and
$\zeta=\sigma_r/\sqrt{2}$, we obtain an optimally conditioned system
\eqref{eq:inc1}, to which the approach of this section can be applied.

\section{Conclusion} \label{AsyRK:sec:conclusion}

We have proposed a simple asynchronous parallel randomized Kaczmarz
algorithm, and proved linear convergence. Our analysis also indicates
the proposed method can be expected to yield near-linear speedup if
the number of processors is bounded by a multiple of the number of
equations in the system. Computational results, including comparison
with an asynchronous stochastic coordinate descent method, confirm the
effectiveness of the approach.

\section*{Acknowledgements}
The authors acknowledge support of National Science Foundation Grants
DMS-0914524 and DMS-1216318, ONR Award N00014-13-1-0129, AFOSR Award
FA9550-13-1-0138, and a Wisconsin Alumni Research Foundation 2011-12
Fall Competition Award. The authors would like to sincerely thank
Yijun Huang for her implementation of algorithms $\AsyRK$ and
$\AsySCD$ used in this paper.

\appendix
\section{Proofs} \label{AsyRK:sec:proofs}

This section provides proofs for our main results in
Section~\ref{sec:mainresults}.

In Algorithm~\ref{alg_rka}, the indices $i(0), t(0), i(1), t(1),
\dotsc, i(j), t(j), \dotsc$ are random variables. We denote the
expectation over all random variables as $\E$, the conditional
expectation with respect to $i(j)$ given $i(0), t(0), i(1), t(1),
\dotsc, i(j-1), t(j-1)$ as $\E_{i(j)}$ and the conditional expectation
with respect to $t(j)$ given $i(0), t(0), i(1), t(1), \dotsc, i(j-1),
t(j-1), i(j)$ as $\E_{t(j)}$.

\medskip

\subsection*{\bf Proof of Theorem~\ref{AsyRK:thm_2}}
\begin{proof}
We start with the following useful results, noting that the random
variable $t(j)$ is distributed uniformly over the set $\mbox{\rm supp}
(a_{i(j)})$:
\begin{equation} \label{eq:useful.1}
\E_{t(j)}(P_{t(j)}a_{i(j)}) ={1\over \theta_{i(j)}}P_{{\rm
    supp}(a_{i(j)})}a_{i(j)} = {1\over \theta_{i(j)}} a_{i(j)}.
\end{equation}

We prove each of the two inequalities in \eqref{AsyRK_eqn_thm2_1} by
induction. We start from the right-hand inequality. First we consider
the expansion of $\|Ax_{j+1}-b\|^2$ for any values of $j$:
\begin{align}
\|Ax_{j+1}-b\|^2 & = \|Ax_j-b\|^2 + \gamma^2\|A\theta_{i(j)}P_{t(j)}a_{i(j)}(a_{i(j)}^Tx_{k(j)}-b_{i(j)})\|^2- \nonumber\\
&\quad\quad 2\gamma\langle Ax_j-b, A\theta_{i(j)}P_{t(j)}a_{i(j)}(a_{i(j)}^Tx_{k(j)}-b_{i(j)}) \rangle \nonumber\\
&=\|Ax_j-b\|^2 + \gamma^2\underbrace{\|A\theta_{i(j)}P_{t(j)}a_{i(j)}(a_{i(j)}^Tx_{k(j)}-b_{i(j)})\|^2}_{T_1} - \nonumber\\
&\quad\quad 2\gamma\underbrace{\langle Ax_{k(j)}-b, A\theta_{i(j)}P_{t(j)}a_{i(j)}(a_{i(j)}^Tx_{k(j)}-b_{i(j)})\rangle}_{T_2} + \nonumber\\
&\quad\quad 2\gamma\underbrace{\langle A(x_{k(j)}-x_{j}), A\theta_{i(j)}P_{t(j)}a_{i(j)}(a_{i(j)}^Tx_{k(j)}-b_{i(j)})\rangle}_{T_3}.
\label{AsyRK:eqn_proof2_1}
\end{align}
Next we consider the expectation of three terms $T_1$, $T_2$, and $T_3$ in
\eqref{AsyRK:eqn_proof2_1}. For $T_1$, we have
\begin{align}
\nonumber
\E(T_1) & =\E(\|A\theta_{i(j)}P_{t(j)}a_{i(j)}(a_{i(j)}^Tx_{k(j)}-b_{i(j)})\|^2) \\
\nonumber
& \leq \alpha^2\E(\|a_{i(j)}^Tx_{k(j)}-b\|^2) \\
\label{eqn_proof2_E1}
& = {\alpha^2 \over m}\E(\|Ax_{k(j)}-b\|^2).
\end{align}
For $T_2$, we have
\begin{align}
\E(T_2) &=\E \langle Ax_{k(j)}-b, A\theta_{i(j)}P_{t(j)}a_{i(j)}(a_{i(j)}^Tx_{k(j)}-b_{i(j)})\rangle
\nonumber \\
&=\E \langle Ax_{k(j)}-b, A \E_{t(j)} (\theta_{i(j)}P_{t(j)}a_{i(j)})(a_{i(j)}^Tx_{k(j)}-b_{i(j)})\rangle
\nonumber \\
&=\E\langle Ax_{k(j)}-b, Aa_{i(j)}(a_{i(j)}^Tx_{k(j)}-b_{i(j)})\rangle
\nonumber \quad \mbox{(by \eqref{eq:useful.1})}
\\
&=\E\langle A^T(Ax_{k(j)}-b), a_{i(j)}(a_{i(j)}^Tx_{k(j)}-b_{i(j)})\rangle
\nonumber \\
&={1\over m}\E(\|A^T(Ax_{k(j)}-b)\|^2).
\label{eqn_proof2_E2}
\end{align}
For $T_3$, we have
\begin{align}
&\E(T_3)=\E\langle A(x_{k(j)}-x_{j}), A\theta_{i(j)}P_{t(j)}a_{i(j)}(a_{i(j)}^Tx_{k(j)}-b_{i(j)})\rangle
\nonumber \\
&\quad =\gamma\E\langle A\sum_{d=k(j)}^{j-1}\theta_{i(d)}P_{t(d)}a_{i(d)}(a^T_{i(d)}x_{k(d)}-b_{i(d)}), A\theta_{i(j)}P_{t(j)}a_{i(j)}(a_{i(j)}^Tx_{k(j)}-b_{i(j)})\rangle
\nonumber \\
&\quad =\gamma\sum_{d=k(j)}^{j-1}\E\langle A\theta_{i(d)}\E_{t(d)}(P_{t(d)})a_{i(d)}(a^T_{i(d)}x_{k(d)}-b_{i(d)}), A\theta_{i(j)}\E_{t(j)}(P_{t(j)})a_{i(j)}(a_{i(j)}^Tx_{k(j)}-b_{i(j)})\rangle
\nonumber \\
&\quad =\gamma\sum_{d=k(j)}^{j-1}\E\langle Aa_{i(d)}(a^T_{i(d)}x_{k(d)}-b_{i(d)}), Aa_{i(j)}(a_{i(j)}^Tx_{k(j)}-b_{i(j)})\rangle
\nonumber \\
&\quad =\gamma\sum_{d=k(j)}^{j-1}\E\langle A\E_{i(d)}(a_{i(d)}(a^T_{i(d)}x_{k(d)}-b_{i(d)})), A\E_{i(j)}(a_{i(j)}(a_{i(j)}^Tx_{k(j)}-b_{i(j)}))\rangle
\nonumber \\
&\quad ={\gamma\over m^2}\sum_{d=k(j)}^{j-1}\E\langle AA^T(Ax_{k(d)}-b), AA^T(Ax_{k(j)}-b)\rangle
\nonumber \\
&\quad \leq{\gamma\over 2m^2}\sum_{d=k(j)}^{j-1}\E(\|AA^T(Ax_{k(d)}-b)\|^2 + \|AA^T(Ax_{k(j)}-b)\|^2)
\nonumber \\
&\quad \leq{\gamma\Lmax^2\over 2m^2}\sum_{d=k(j)}^{j-1} 
\left[ \E(\|Ax_{k(d)}-b\|^2 + \|Ax_{k(j)}-b\|^2) \right],
\label{eqn_proof2_E3}
\end{align}
where the third line is from the observation that $t(d)$ and $t(j)$
are conditionally independent given $i(d)$ and $i(j)$; the fifth line
uses the result that $i(d)$ only affects $x_{d+1}$ and subsequent
iterates and $k(j)$ is less than $d+1$ (so $x_{k(j)}$ and $i(d)$ are
independent to each other). Combining \eqref{eqn_proof2_E1},
\eqref{eqn_proof2_E2}, \eqref{eqn_proof2_E3}, and
\eqref{AsyRK:eqn_proof2_1}, we obtain
\begin{align}
\E(\|Ax_{j+1}-b\|^2) & \leq \E(\|Ax_j-b\|^2) + {\alpha^2 \gamma^2\over m}\E(\|Ax_{k(j)}-b\|^2) - {2\gamma\over m}\E(\|A^T(Ax_{k(j)}-b)\|^2) +
\nonumber \\
& \quad\quad {\gamma^2\Lmax^2\over m^2}\sum_{d=k(j)}^{j-1} \left[ \E(\|Ax_{k(d)}-b\|^2 + \|Ax_{k(j)}-b\|^2) \right].
\label{eqn_proof2_obj}
\end{align}
We can use this bound to show that the right-hand inequality
in~\eqref{AsyRK_eqn_thm2_1} holds for $j=0$. By setting $j=0$
in~\eqref{eqn_proof2_obj} and noting that $k(0)=0$ and that the last
summation is vacuous, we obtain
\begin{align}
\E(\|Ax_{1}-b\|^2) & \leq \E(\|Ax_0-b\|^2) + {\gamma^2\alpha^2 \over m}\E(\|Ax_{0}-b\|^2) - {2\gamma\over m}\E(\|A^T(Ax_{0}-b)\|^2) 
\nonumber \\
&\leq \|Ax_0-b\|^2 + {\gamma^2\alpha^2 \over m} \|Ax_{0}-b\|^2
\nonumber \\
&\leq \left(1+{\gamma^2\alpha^2 \over m}\right) \|Ax_0-b\|^2.
\label{eqn_proof2_j0}
\end{align}
From the third bound in~\eqref{eqn_thm2_gamma}, we have
\[
1+{\gamma^2\alpha^2 \over m} \leq 1 + \frac{m(\rho-1)\alpha^2}{\rho^{\tau}(m\alpha^2+\Lmax^2\tau\rho^\tau)} \leq 1 + \frac{m(\rho-1)\alpha^2}{\rho^{\tau}m\alpha^2} \leq 1 + (\rho-1) = \rho,
\]
where the third inequality follows from $\rho>1$. By substituting
into~\eqref{eqn_proof2_j0}, we obtain $\E(\|x_0-x_0^*\|^2)\leq \rho
\E(\|x_1-x_1^*\|^2)$.

For the inductive step, we use~\eqref{eqn_proof2_obj} again, assuming
that the right-hand inequality in~\eqref{AsyRK_eqn_thm2_1} holds up to
stage $j$, and thus that
\[
\E(\|Ax_{k(j)}-b\|^2) \leq \rho^{\tau}\E(\|Ax_{j}-b\|^2)\quad\text{and}\quad\E(\|Ax_{k(d)}-b\|^2) \leq \rho^{2\tau}\E(\|Ax_{j}-b\|^2)
\]
provided that $0\leq j-k(j)\leq \tau$ and $0\leq j-k(d)\leq 2\tau$, as
assumed. By substituting into the right-hand side
of~\eqref{eqn_proof2_obj} again, we obtain
\begin{align}
\E(\|Ax_{j+1}-b\|^2) &\leq \E(\|Ax_j-b\|^2) + {\alpha^2 \gamma^2\over m}\E(\|Ax_{k(j)}-b\|^2)  +
\nonumber \\
&\quad\quad{\gamma^2\Lmax^2\over m^2}\sum_{d=k(j)}^{j-1}\E(\|Ax_{k(d)}-b\|^2 + \|Ax_{k(j)}-b\|^2)
\nonumber \\
&\leq \E(\|Ax_j-b\|^2) + {\rho^{\tau}\alpha^2 \gamma^2\over m}\E(\|Ax_{j}-b\|^2)  +
\nonumber \\
&\quad\quad{\gamma^2\Lmax^2\over m^2}\sum_{d=k(j)}^{j-1}\left(\rho^{2\tau}\E(\|Ax_{j}-b\|^2) + \rho^{\tau}\E(\|Ax_{j}-b\|^2)\right)
\nonumber \\
&\leq \E(\|Ax_j-b\|^2) + {\rho^{\tau}\alpha^2 \gamma^2\over m}\E(\|Ax_{j}-b\|^2)  +{\gamma^2\Lmax^2\over m^2}(2 \tau \rho^{2\tau}) \E(\|Ax_{j}-b\|^2)
\nonumber \\
&= \left(1+\gamma^2\left({\rho^{\tau}(m\alpha^2 +2\tau\Lmax^2\rho^{\tau})\over m^2}\right)\right)\E(\|Ax_j-b\|^2)
\nonumber \\
&\leq (1+(\rho-1))\E(\|Ax_j-b\|^2)
\nonumber \\
&= \rho\E(\|Ax_j-b\|^2),
\nonumber
\end{align}
where the last inequality uses the third bound on $\gamma$
from~\eqref{eqn_thm2_gamma}. We conclude that the right-hand side
inequality in~\eqref{AsyRK_eqn_thm2_1} holds for all $j$.

We now work on the left-hand inequality
in~\eqref{AsyRK_eqn_thm2_1}. For all $j$, we have the following:
%
\begin{align}
&\E(\|Ax_{j+1}-b\|^2) = \E(\|Ax_j-b\|^2) + \gamma^2\E(\|A\theta_{i(j)}P_{t(j)}a_{i(j)}(a_{i(j)}^Tx_{k(j)}-b_{i(j)})\|^2)-
\nonumber\\
&\quad\quad
2\gamma\E(\langle Ax_j-b, A\theta_{i(j)}P_{t(j)}a_{i(j)}(a_{i(j)}^Tx_{k(j)}-b_{i(j)}) \rangle)
\nonumber\\
&\quad
\geq \E(\|Ax_j-b\|^2) - 2\gamma\E(\langle Ax_j-b, A\theta_{i(j)}P_{t(j)}a_{i(j)}(a_{i(j)}^Tx_{k(j)}-b_{i(j)}) \rangle)
\nonumber\\
&\quad
\geq \E(\|Ax_j-b\|^2) - 2\gamma\E(\langle Ax_j-b, A \E_{t(j)} (\theta_{i(j)}P_{t(j)}a_{i(j)})(a_{i(j)}^Tx_{k(j)}-b_{i(j)}) \rangle)
\nonumber\\
&\quad
= \E(\|Ax_j-b\|^2) - 2\gamma\E\langle A^T(Ax_{j}-b), a_{i(j)}(a_{i(j)}^Tx_{k(j)}-b_{i(j)})\rangle \quad\quad \mbox{(using \eqref{eq:useful.1})}
\nonumber\\
&\quad
= \E(\|Ax_j-b\|^2) - {2\gamma\over m}\E\langle A^T(Ax_{j}-b), A^T(Ax_{k(j)}-b)\rangle
\nonumber \\
&\quad
\geq E(\|Ax_j-b\|^2) - {\gamma\over m}\E( \|A^T(Ax_{j}-b)\|^2 + \|A^T(Ax_{k(j)}-b)\|^2)
\nonumber \\
&\quad
\geq E(\|Ax_j-b\|^2) - {\gamma\Lmax\over m}\E( \|Ax_{j}-b\|^2 + \|(Ax_{k(j)}-b\|^2)
\nonumber \\
&\quad
= \left(1 - {\gamma\Lmax\over m}\right)\E( \|Ax_{j}-b\|^2) - {\gamma\Lmax\over m} \E(\|Ax_{k(j)}-b\|^2).
\label{eqn_proof2_obj2}
\end{align}
We can use this bound to show that the left-hand inequality
in~\eqref{AsyRK_eqn_thm2_1} holds for $j=0$. By setting $j=0$
in~\eqref{eqn_proof2_obj} and noting that $k(0)=0$, we obtain
\[
\E(\|Ax_{1}-b\|^2) \geq \left(1 - {2\gamma\Lmax\over m}\right)\E(\|Ax_{0}-b\|^2).
\]
From the second bound in~\eqref{eqn_thm2_gamma}, we have
\[
1 - {2\gamma\Lmax\over m} \geq 1 - {\rho-1 \over \rho^{\tau+1}} = 1 - {1-\rho^{-1} \over \rho^{\tau}} \geq \rho^{-1},
\]
where the last inequality follows from $\rho > 1$. By substituting
into~\eqref{eqn_proof2_obj2}, we obtain $\rho^{-1}\E(\|Ax_0-b\|^2)\leq
\E(\|Ax_1-b\|^2)$. For the inductive step, we
use~\eqref{eqn_proof2_obj2} again, assuming that the left-hand
inequality in~\eqref{AsyRK_eqn_thm2_1} holds up to stage $j$, and thus
that
\[
\E(\|Ax_j-b\|^2) \geq \rho^{-\tau}\E(\|Ax_{k(j)}-b\|^2),
\]
provided that $0\leq j-k(j)\leq \tau$, as assumed \eqref{eq:kj}.  By
substituting into the left-hand side of~\eqref{eqn_proof2_obj2} again,
we obtain
\begin{align}
\E(\|Ax_{j+1}-b\|^2) & \geq  \left(1 - {\gamma\Lmax\over m}\right)\E( \|Ax_{j}-b\|^2) - {\gamma\rho^{\tau}\Lmax\over m} \E(\|Ax_{j}-b\|^2)
\nonumber \\ 
&\geq  \left(1 - {2\gamma\rho^{\tau}\Lmax\over m}\right)\E( \|Ax_{j}-b\|^2).
\label{eqn_proof2_jj}
\end{align}
From the second bound in~\eqref{eqn_thm2_gamma}, we have
\[
1 - {2\gamma\rho^{\tau}\Lmax\over m} \geq 1 - {\rho-1\over \rho} = \rho^{-1}.
\]
We conclude that the left-hand side inequality in~\eqref{AsyRK_eqn_thm2_1} holds for all $j$.

At this point, we have shown that both inequalities
in~\eqref{AsyRK_eqn_thm2_1} are satisfied for all $j$.


We next prove \eqref{AsyRK:eqn_thm2_2}. Consider the expansion of
$\|x_{j+1}-x_{j+1}^*\|^2$:
\begin{equation}
\begin{aligned}
&\|x_{j+1}-x_{j+1}^*\|^2 = \|x_j-\gamma \theta_{i(j)}P_{t(j)}a_{i(j)}(a_{i(j)}^Tx_{k(j)}-b_{i(j)}) - x_{j+1}^*\|^2\\
&\quad \leq  \|x_j-\gamma \theta_{i(j)}P_{t(j)}a_{i(j)}(a_{i(j)}^Tx_{k(j)}-b_{i(j)}) - x_{j}^*\|^2\\
&\quad =\|x_j-x_j^*\|^2 + \gamma^2\theta_{i(j)}^2\|P_{t(j)}a_{i(j)}(a^T_{i(j)}x_{k(j)}-b_{i(j)})\|^2 - 2\gamma\langle x_j-x_j^*, \theta_{i(j)} P_{t(j)}a_{i(j)}(a_{i(j)}^Tx_{k(j)}-b_{i(j)}) \rangle\\
&\quad =\|x_j-x_j^*\|^2 + \gamma^2\underbrace{\theta_{i(j)}^2\|P_{t(j)}a_{i(j)}(a^T_{i(j)}x_{k(j)}-b_{i(j)})\|^2}_{T_4} - \\
&\quad \quad 2\gamma\underbrace{\langle x_{k(j)} -x_j^*, \theta_{i(j)} P_{t(j)}a_{i(j)}(a_{i(j)}^Tx_{k(j)}-b_{i(j)}) \rangle}_{T_5} +2\gamma \underbrace{\langle x_{k(j)} -x_j, \theta_{i(j)} P_{t(j)}a_{i(j)}(a_{i(j)}^Tx_{k(j)}-b_{i(j)}) \rangle}_{T_6},
\end{aligned}
\label{AsyRK:eqn_proof2_6}
\end{equation}
Next, we estimate the expectations of $T_4$, $T_5$, and $T_6$. For
$T_4$, we have
\begin{align}
\nonumber
\E(T_4)& =\E(\theta_{i(j)}^2\|P_{t(j)}a_{i(j)}(a_{i(j)}^Tx_{k(j)}-b_{i(j)})\|^2)\\
\nonumber
&= \E(\theta_{i(j)}^2\E_{t(j)}(\|P_{t(j)}a_{i(j)}(a_{i(j)}^Tx_{k(j)}-b_{i(j)})\|^2) \\
\nonumber
&= \theta_{i(j)}^2 \E \left( \frac{1}{\theta_{i(j)}} \sum_{t \in \mbox{\rm supp}(a_{i(j)})} 
(a_{i(j)}^T x_{k(j)} - b_{i(j)})^2 a_{i(j)}^T P_t a_{i(j)} \right) \\
\nonumber
&= \E(\theta_{i(j)}\|a_{i(j)}(a_{i(j)}^Tx_{k(j)}-b_{i(j)})\|^2) \\
\nonumber
&\leq \mu \E(a_{i(j)}^Tx_{k(j)}-b_{i(j)})^2\\
\label{eqn_proof2_7}
&= {\mu\over m} \E(\|Ax_{k(j)}-b\|^2),
\end{align}
For $T_5$, we have
\begin{align}
\nonumber
\E(T_5)& =\E(\langle x_{k(j)} -x_j^*, \theta_{i(j)} P_{t(j)}a_{i(j)}(a_{i(j)}^Tx_{k(j)}-b_{i(j)}) \rangle) \\
\nonumber
&= \E(\langle x_{k(j)} -x_j^*, \E_{t(j)} (\theta_{i(j)} P_{t(j)}a_{i(j)}) (a_{i(j)}^Tx_{k(j)}-b_{i(j)}) \rangle) \\
\nonumber
&=\E(\langle x_{k(j)} -x_j^*, a_{i(j)}(a_{i(j)}^Tx_{k(j)}-b_{i(j)}) \rangle)\\
\label{eqn_proof2_8}
&={1\over m}\E(\|Ax_{k(j)}-b\|^2).
\end{align}
By following a derivation similar to \eqref{eqn_proof2_E3} for $T_6$,
we obtain 
\begin{align}
\nonumber
\E(T_6)&=\E(\langle x_{k(j)} -x_j, \theta_{i(j)} P_{t(j)}a_{i(j)}(a_{i(j)}^Tx_{k(j)}-b_{i(j)})\rangle) \\
\nonumber
&=\E\left(\gamma\sum_{d=k(j)}^{j-1}\langle \theta_{i(d)}P_{t(d)}a_{i(d)}(a^T_{i(d)}x_{k(d)}-b_{i(d)}), \theta_{i(j)} P_{t(j)}a_{i(j)}(a_{i(j)}^Tx_{k(j)}-b_{i(j)})\rangle\right) \\
\nonumber
&={\gamma\over m^2}\E\left(\sum_{d=k(j)}^{j-1}\langle A^T(Ax_{k(d)}-b), A^T(Ax_{k(j)}-b) \rangle\right) \\
\nonumber
&\leq{\gamma\over m^2}\E\left(\sum_{d=k(j)}^{j-1}\left|\langle A^T(Ax_{k(d)}-b), A^T(Ax_{k(j)}-b) \rangle\right|\right) \\
\nonumber
&\leq{\gamma\over 2m^2}\E\left(\sum_{d=k(j)}^{j-1}\|A^T(Ax_{k(d)}-b)\|^2+\|A^T(Ax_{k(j)}-b)\|^2\right)\\
\label{eq:crap2}
&\leq{\gamma \Lmax\over 2m^2}\E\left(\sum_{d=k(j)}^{j-1}\|Ax_{k(d)}-b\|^2+\|Ax_{k(j)}-b\|^2\right).
\end{align}
Since for $d=k(j), k(j)+1, \dotsc, j-1$, we have
\[
k(j)-\tau \le k(d) \le j-2 \le k(j)+\tau-1,
\]
it follows from \eqref{AsyRK_eqn_thm2_1} that
\[
\| A x_{k(d)} - b \|^2 \le \rho^{\tau} \| A x_{k(j)}-b \|^2)
\]
Thus from \eqref{eq:crap2}, we have
\begin{equation}
\label{eqn_proof2_10}
\E(T_6) \leq{\gamma \tau \Lmax (1+\rho^{\tau})\over 2m^2}\E(\|Ax_{k(j)}-b\|^2)
\leq {\gamma\tau\Lmax\rho^{\tau}\over m^2}\E(\|Ax_{k(j)}-b\|^2).
\end{equation}
By substituting \eqref{eqn_proof2_7}, \eqref{eqn_proof2_8}, and
\eqref{eqn_proof2_10} into \eqref{AsyRK:eqn_proof2_6}, we obtain
\[
\E(\|x_{j+1}-x_{j+1}^*\|^2 ) \leq \E(\|x_j-x_j^*\|^2) - \left({2\gamma
  - \mu\gamma^2\over m}-{2\gamma^2\tau\Lmax\rho^{\tau}\over
  m^2}\right)\E(\|Ax_{k(j)}-b\|^2).
\]
Since by \eqref{eqn_thm_psi} and \eqref{eqn_thm2_gamma}, we have
\[
{2\gamma - \mu\gamma^2\over m}-{2\gamma^2\tau\Lmax\rho^{\tau}\over
  m^2} =
\frac{\gamma}{m} \left[ 2- \gamma \left(\mu + \frac{2 \tau \Lmax \rho^{\tau}}{m}\right)\right] = \frac{\gamma}{m} (2-\gamma \psi) >0,
\]
we have  from the bound above that
\begin{align*}
\E(\|x_{j+1}-x_{j+1}^*\|^2 ) &
\leq \E(\|x_j-x_j^*\|^2) - \frac{\gamma}{m} (2-\gamma \psi)  \E(\|Ax_{k(j)}-b\|^2)\\
&\leq \E(\|x_j-x_j^*\|^2) - {\gamma\over m}(2-\psi\gamma) \Lmin \E(\|x_{k(j)}-x_{k(j)}^*\|^2)\\
&\leq \left(1-{\Lmin\gamma\over m}(2-\psi\gamma)\right)\E(\|x_j-x_j^*\|^2),
\end{align*}
where the second line implies that $\E(\|x_j-x_j^*\|^2)$ is
monotonically decreasing, and the last line is obtained by the
implication from the second line. This completes the proof of
\eqref{AsyRK:eqn_thm2_2}.
\end{proof}

\noindent{\bf Proof of Corollary~\ref{cor}}
\begin{proof}
Note first that for $\rho$ defined by~\eqref{eqn_cor_rho}, and using
\eqref{eqn_cor_ass}, we have
\begin{align*}
\rho^{\tau} \leq \rho^{\tau+1} = \left[\left(1+{2e\Lmax\over
    {m}}\right)^{m\over 2e\Lmax} \right]^{{2e\Lmax(\tau+1)\over {m}}}
\leq e^{{2e\Lmax(\tau+1)\over {m}}} \leq e.
\end{align*}
Thus from the definition of $\psi$~\eqref{eqn_thm_psi}, and using
\eqref{eqn_cor_ass} again, we have
\begin{align}
\psi = \mu + \frac{2\Lmax \tau\rho^{\tau}}{m} \leq \mu + \frac{2e\tau\Lmax}{m} \leq \mu +1.
\label{eqn_cor_proof_psibound}
\end{align}
We show now that the steplength parameter choice $\gamma={1/\psi}$
satisfies all the bounds in~\eqref{eqn_thm2_gamma}, by showing that
the second and third bounds are implied by the first. For the second
bound, we have
\[
\frac{m(\rho-1)}{2\Lmax\rho^{\tau+1}} \geq \frac{m(\rho-1)}{2e\Lmax}
\geq \tau + 1 \geq 1 \ge \frac{1}{\psi},
\]
where the second inequality follows from~\eqref{eqn_cor_rho} and the
final inequality follows from the definition of $\psi$ in
\eqref{eqn_thm_psi} and the fact that $\psi > \mu \ge 1$.

For the third bound in~\eqref{eqn_thm2_gamma}, we have (by taking
squares) that
\begin{alignat*}{2}
&m^2 \frac{(\rho-1)}{\rho^{\tau}(m\alpha^2+\Lmax^2\tau\rho^{\tau})} 
= \frac{2me\Lmax (\tau+1)}{\rho^{\tau}(m\alpha^2+\Lmax^2\tau\rho^{\tau})} 
&&\quad \left(\text{from the definition of $\rho$ in \eqref{eqn_cor_rho}}\right)\\
&\quad \geq \frac{2me\Lmax (\tau+1)}{e(m\Lmax\mu^2+\Lmax^2\tau e)} \\
&\quad = \frac{1}{\frac{\mu^2}{2(\tau+1)} + \frac{\Lmax \tau e}{2m(\tau+1)}}\\
&\quad \geq \frac{1}{{\frac{\mu^2}{2(\tau+1)}} + \frac{\tau}{4(\tau+1)^2}} 
&&\quad \left(\text{from the lower bound of $m$ in \eqref{eqn_cor_ass}}\right)\\
&\quad \geq \frac{1}{\frac{\mu^2}{2} + \frac{1}{16}} 
&& \quad \left(\text{from $\frac{\tau}{(\tau+1)^2}\leq \frac14$}\right) \\
&\quad \geq \frac{1}{\mu^2} \ge \frac{1}{\psi^2}.
\end{alignat*}
We can thus set $\gamma=1/\psi$, and by substituting this choice
into~\eqref{AsyRK:eqn_thm2_2} and
using~\eqref{eqn_cor_proof_psibound}, we
obtain~\eqref{eqn_cor_estbound}.
\end{proof}

\noindent{\bf Proof of Theorem~\ref{thm_3}}
\begin{proof}
From Markov's inequality, we have
\begin{align*}
\mathbb{P}(\|x_{j}-x_j^*\|^2\geq \epsilon) & \leq \epsilon^{-1}\mathbb{E}(\|x_{j}-x_j^*\|^2) \\
& \leq \epsilon^{-1}\left(1-{\Lmin\over m(\mu+1)}\right)^{j}\|x_0-x_0^*\|^2\\
&\le
{\epsilon}^{-1}{\left(1-c\right)^{(1/c) \left|{\log{\|x_0-x_0^*\|^2\over \eta\epsilon}}\right|}\|x_0-x_0^*\|^2} \quad \left(\mbox{with $c={\Lmin\over m(\mu+1)}$}\right) \\
&\le
{\epsilon}^{-1}\|x_0-x_0^*\|^2 e^{-\left|{\log{\|x_0-x_0^*\|^2\over \eta\epsilon}}\right|} \\
&=
{\eta}e^{\log{\|x_0-x_0^*\|^2\over \eta\epsilon}} e^{-\left|{\log{\|x_0-x_0^*\|^2\over \eta\epsilon}}\right|}\\
&\le
\eta,
\end{align*}
where the second inequality applies~\eqref{eqn_cor_estbound}, the third inequality uses the definition of $j$~\eqref{eqn_thm3_j}, and the second last inequality uses the inequality $(1-c)^{1/c}\leq e^{-1}~\forall c\in (0, 1)$, which completes the proof.
\end{proof}

{
\bibliographystyle{plainnat}
\bibliography{reference}

\begin{thebibliography}{27}
\providecommand{\natexlab}[1]{#1}
\providecommand{\url}[1]{\texttt{#1}}
\expandafter\ifx\csname urlstyle\endcsname\relax
  \providecommand{\doi}[1]{doi: #1}\else
  \providecommand{\doi}{doi: \begingroup \urlstyle{rm}\Url}\fi

\bibitem[Aharoni and Censor(1989)]{Aharoni89}
R.~Aharoni and Y.~Censor.
\newblock Block-iterative projection methods for parallel computation of
  solutions to convex feasibility problems.
\newblock \emph{Linear Algebra and Its Applications}, 120:\penalty0 165--175,
  1989.

\bibitem[{Avron} et~al.(2014){Avron}, {Druinsky}, and {Gupta}]{Avron13arXiv}
H.~{Avron}, A.~{Druinsky}, and A.~{Gupta}.
\newblock Revisiting asynchronous linear solvers: Provable convergence rate
  through randomization.
\newblock \emph{IPDPS}, 2014.

\bibitem[Bertsekas and Tsitsiklis(1989)]{Bertsekas89}
D.~P. Bertsekas and J.~N. Tsitsiklis.
\newblock \emph{Parallel and distributed computation: numerical methods}.
\newblock Pentice Hall, 1989.

\bibitem[{C}ensor et~al.(2001){C}ensor, {G}ordon, and {G}ordon]{Cegogo01}
{Y.} {C}ensor, {D.} {G}ordon, and {R.} {G}ordon.
\newblock {C}omponent averaging: {A}n efficient iterative parallel algorithm
  for large and sparse unstructured problems.
\newblock \emph{{P}arallel {C}omputing}, 27\penalty0 (6):\penalty0 777--808,
  2001.

\bibitem[Eggermont(1981)]{Eggermont81}
P.~P.~B. Eggermont.
\newblock Iterative algorithms for large partitioned linear systems, with
  applications to image reconstruction.
\newblock \emph{Linear Algebra and its Applications}, 40:\penalty0 37--67,
  1981.

\bibitem[Eldar and Needell(2011)]{EldarN11}
Y.~C. Eldar and D.~Needell.
\newblock Acceleration of randomized {Kaczmarz} method via the
  {Johnson-Lindenstrauss} lemma.
\newblock \emph{Numerical Algorithms}, 58\penalty0 (2):\penalty0 163--177,
  2011.

\bibitem[Elfving and Nikazad(2009)]{ElfvingNikazad09}
T.~Elfving and T.~Nikazad.
\newblock Properties of a class of block-iterative methods.
\newblock \emph{Inverse problems}, 25, 2009.

\bibitem[Elsner et~al.(1990)Elsner, Koltracht, and Neumann]{Elsner90}
L.~Elsner, I.~Koltracht, and M.~Neumann.
\newblock On the convergence of asynchronous paracontractions with application
  to tomographic reconstruction from incomplete data.
\newblock \emph{Linear Algebra and its Applications}, 130:\penalty0 65--82,
  1990.

\bibitem[Ferris and Mangasarian(1994)]{Ferris94}
M.~C. Ferris and O.~L. Mangasarian.
\newblock Parallel variable distribution.
\newblock \emph{SIAM Journal on Optimization}, 4\penalty0 (4):\penalty0
  815--832, 1994.

\bibitem[Galantai(2005)]{Galantai25}
A.~Galantai.
\newblock On the rate of convergence of the alternating projection method in
  finite dimensional spaces.
\newblock \emph{Journal of Mathematical Analysis and Applications},
  310:\penalty0 30--44, 2005.

\bibitem[Herman(1980)]{Herman80}
G.~T. Herman.
\newblock \emph{Image Reconstruction from Projections: The Fundamentals of
  Computerized Tomography}.
\newblock Academic Press, 1980.

\bibitem[Herman(2009)]{Herman09}
G.~T. Herman.
\newblock \emph{Fundamentals of Computerized Tomography}.
\newblock Springer, 2009.

\bibitem[Hoffman(1952)]{Hoffman52}
A.~J. Hoffman.
\newblock On approximate solutions of systems of linear inequalities.
\newblock \emph{Journal of Research of the National Bureau of Standards},
  49\penalty0 (4):\penalty0 263--265, 1952.

\bibitem[Kaczmarz(1937)]{Kaczmarz37}
S.~Kaczmarz.
\newblock Angenaherte auflsung von systemen linearer gleichungen.
\newblock \emph{Bulletin International de l'Acadmie Polonaise des Sciences et
  des Letters}, 35:\penalty0 355--357, 1937.

\bibitem[Leventhal and Lewis(2010)]{Leventhal08}
D.~Leventhal and A.~S. Lewis.
\newblock Randomized methods for linear constraints: Convergence rates and
  conditioning.
\newblock \emph{Mathematics of Operations Research}, 35\penalty0 (3):\penalty0
  641--654, 2010.

\bibitem[Liu and Wright(2013)]{LiuWright13b}
J.~Liu and S.~J. Wright.
\newblock An accelerated randomized {K}aczmarz algorithm.
\newblock \emph{arXiv:1310.2887}, 2013.

\bibitem[{Liu} et~al.(2013){Liu}, {Wright}, {R{\'e}}, {Bittorf}, and
  Sridhar]{LiuWright13}
J.~{Liu}, S.~J. {Wright}, C.~{R{\'e}}, V.~{Bittorf}, and S.~Sridhar.
\newblock An asynchronous parallel stochastic coordinate descent algorithm.
\newblock \emph{arXiv: 1311.1873}, 2013.

\bibitem[Needell(2010)]{Needell10}
D.~Needell.
\newblock Randomized {Kaczmarz} solver for noisy linear systems.
\newblock \emph{{BIT} Numerical Mathematics}, 50\penalty0 (2):\penalty0
  1422--1436, 2010.

\bibitem[Nemirovski et~al.(2009)Nemirovski, Juditsky, Lan, and
  Shapiro]{Nemirovski09}
A.~Nemirovski, A.~Juditsky, G.~Lan, and A.~Shapiro.
\newblock Robust stochastic approximation approach to stochastic programming.
\newblock \emph{SIAM Journal on Optimization}, 19:\penalty0 1574--1609, 2009.

\bibitem[Nikazad(2008)]{Nikazad08}
T.~Nikazad.
\newblock Algebraic reconstruction methods.
\newblock \emph{Thesis}, 2008.

\bibitem[Niu et~al.(2011)Niu, Recht, R{\'e}, and Wright]{Hogwild11nips}
F.~Niu, B.~Recht, C.~R{\'e}, and S.~J. Wright.
\newblock Hogwild: A lock-free approach to parallelizing stochastic gradient
  descent.
\newblock \emph{NIPS}, pages 693--701, 2011.

\bibitem[Popa(1999)]{Popa99}
C.~Popa.
\newblock Characterization of the solutions set of least-squares problems by an
  extension of {Kaczmarz's} projections method.
\newblock \emph{Journal of Applied Mathematics and Computing}, 6:\penalty0
  51--64, 1999.

\bibitem[{Richt{\'a}rik} and {Tak{\'a}{\v c}}(2012)]{Richtarik12arXiv}
P.~{Richt{\'a}rik} and M.~{Tak{\'a}{\v c}}.
\newblock Parallel coordinate descent methods for big data optimization.
\newblock \emph{ArXiv: 1212.0873}, 2012.

\bibitem[{Sridhar} et~al.(2013){Sridhar}, {Bittorf}, {Liu}, {Zhang}, {R{\'e}},
  and {Wright}]{Sridhar2013nips}
S.~{Sridhar}, V.~{Bittorf}, J.~{Liu}, C.~{Zhang}, C.~{R{\'e}}, and S.~J.
  {Wright}.
\newblock An approximate, efficient solver for {LP} rounding.
\newblock \emph{NIPS}, 2013.

\bibitem[Strohmer and Vershynin(2009)]{Strohmer09}
T.~Strohmer and R.~Vershynin.
\newblock A randomized kaczmarz algorithm with exponential convergence.
\newblock \emph{Journal of Fourier Analysis and Applications}, 15:\penalty0
  262--278, 2009.

\bibitem[Vershynin(2011)]{Vershynin2011arXiv}
R.~Vershynin.
\newblock Introduction to the non-asymptotic analysis of random matrices.
\newblock \emph{ArXiv: 1011.3027}, 2011.

\bibitem[Zouzias and Freris(2012)]{Zouzias12}
A.~Zouzias and N.~M. Freris.
\newblock Randomized extended {Kaczmarz} for solving least-squares.
\newblock \emph{arXiv:1205.5770v2}, 2012.

\end{thebibliography}
}
\end{document}